\documentclass[12pt]{amsart}
\usepackage{a4wide} 
\usepackage{amsmath} 
\usepackage{amsfonts}
\usepackage{amsthm} 
\usepackage{amssymb} 

\theoremstyle{plain}
\newtheorem{thm}{Theorem}
\newtheorem{prop}[thm]{Proposition} 
\newtheorem{lem}[thm]{Lemma} 
\newtheorem{cor}[thm]{Corollary}

\theoremstyle{remark}

\newtheorem*{example*}{Example}

\newtheorem*{rem*}{Remark}

\newcommand{\Z}{\mathbb{Z}}
\newcommand{\Q}{\mathbb{Q}}
\newcommand{\R}{\mathbb{R}}
\newcommand{\C}{\mathbb{C}}

\renewcommand{\H}{\mathbb{H}}
\newcommand{\F}{\mathbb{F}}

\newcommand{\zxz}[4]{\begin{pmatrix} #1 & #2 \\ #3 & #4 \end{pmatrix}}

\newcommand{\kzxz}[4]{\left(\begin{smallmatrix} #1 & #2 \\ #3 & #4\end{smallmatrix}\right) }
\newcommand{\kabcd}{\kzxz{a}{b}{c}{d}}

\newcommand{\calA}{\mathcal{A}}

\newcommand{\calM}{\mathcal{M}}

\newcommand{\calO}{\mathcal{O}}

\newcommand{\calS}{\mathcal{S}}

\newcommand{\frakd}{\mathfrak d}
\newcommand{\frake}{\mathfrak e}

\newcommand{\frakE}{\mathfrak E}

\newcommand{\eps}{\varepsilon}

\newcommand{\norm}{\operatorname{N}}

\newcommand{\tr}{\operatorname{tr}}

\newcommand{\Sl}{\operatorname{SL}}
\newcommand{\Gl}{\operatorname{GL}}

\newcommand{\Orth}{\operatorname{O}}

\newcommand{\Gal}{\operatorname{Gal}}

\begin{document}

\title[On Borcherds products]{On Borcherds products associated with lattices of prime discriminant}

\author[Jan H.~Bruinier and M.~Bundschuh]{Jan Hendrik Bruinier and Michael Bundschuh}
\date{May 29, 2001, revision: December 3, 2001}
\dedicatory{Dedicated to the memory of Robert A. Rankin}
\address{University of Wisconsin-Madison, Department of Mathematics, Van Vleck Hall, 480 Lincoln Drive, Madison, WI 53706-1388, USA}
\email{bruinier@math.wisc.edu} 
\address{Mathematisches Institut, Universtit\"at Heidelberg, Im Neuenheimer Feld 288, D-69120 Heidelberg, Germany}
\email{bundschu@mathi.uni-heidelberg.de }
\thanks{The first author thanks the Number Theory Foundation for their generous support.}
\maketitle

\section{Introduction}

Let $(V,q)$ be a real quadratic space and $L\subset V$ an even lattice with dual $L'$. 
In the description of the theta lifting from elliptic modular forms on $\Sl_2(\R)$ to modular forms for the orthogonal group of $L$ it is often convenient to work with vector valued modular forms for the full group $\Sl_2(\Z)$ rather than with scalar valued modular forms for congruence subgroups. 
In particular, such vector valued modular forms occur naturally in the theory of automorphic products due to Borcherds \cite{Bo1}. Their transformation behavior under $\Sl_2(\Z)$ is dictated by the Weil representation 
associated with the discriminant group $L'/L$.

However, there is no smooth structure theory for vector valued modular forms. For instance, 
they do not form a weight-graded algebra, and there is no natural action of the full Hecke algebra and the Galois group over $\Q$.
 
In the present note, we 
consider the special case that the discriminant group has odd prime order $p$. 
This implies that the lattice has even dimension.
We show that the relevant spaces of vector valued modular forms can be described using scalar valued modular forms for the group $\Gamma_0(p)$, whose Fourier expansion is supported on either the squares or the non-squares modulo $p$ (see Theorem \ref{mainres}).
In the proof we use some basic properties of the Weil representation and an idea due to Krieg, who considered the special case that $L$ is given by the ring of integers in an imaginary quadratic field \cite{Kr}. Related results for certain lattices of odd dimension were obtained earlier by Eichler, Zagier and Skoruppa (cf.~\cite{EZ} chapter 5, \cite{Sk}).

If the signature of $L$ is $(2,n)$ one obtains a very explicit description of the Borcherds lifting: It maps scalar valued modular forms of weight $1-n/2$ as above which are holomorphic except for a pole at the cusp $\infty$ to meromorphic modular forms for the orthogonal group attached to $L$. 
Furthermore, using the Serre-duality result of \cite{Bo2} one gets an existence criterion for Borcherds products involving a similar space of cusp forms  of weight $1+n/2$ for $\Gamma_0(p)$ (see Theorem \ref{crit}). In section \ref{sec4} we present these ideas in the special $\Orth(2,2)$-case of Hilbert modular surfaces of prime discriminant (Theorem \ref{hilbert}). Since we feel that this case is of particular interest and since there are  no explicit examples in the literature so far (to the best of our knowledge), we work out some Borcherds products in detail.  In the computation of the Weyl vectors we need to apply a result of \cite{Br2}.
For instance, we construct the product of theta-series for $\Q(\sqrt{5})$ considered by Gundlach in \cite{Gu} as a Borcherds product.

\section{Vector valued modular forms}\label{sec2}

We begin by fixing some notation.
We denote by $L$ an even lattice of signature $(b^+,b^-)$ equipped with a non-degenerate quadratic form $q(x)=\frac{1}{2}(x,x)$. We write $L'$ for the dual lattice, $m=b^++b^-$ for the dimension, and put $r=b^+-b^-$.

Throughout we assume that the discriminant group $L'/L$ has prime order $p\geq 3$.
This implies that $m$ and $r$ are even.
If $S$ denotes the Gram matrix of $L$, then it is well known that
$(-1)^{m/2} \det(S)= (-1)^{r/2} p$
is a discriminant and therefore congruent to $1$ modulo $4$. Thus $p$ determines $r$  modulo $4$. 

The following argument shows that $p$ modulo $4$ together with the type of the quadratic form on $L'/L$ induced by $q$ determine $r$ modulo $8$:
Let $\chi_p(\cdot)=(\frac{\cdot}{p})$ denote the Dirichlet character given by the Legendre symbol and define
\[
\eps_p=\begin{cases}1,& p\equiv 1 \pmod{4},\\
i,& p\equiv 3\pmod{4}.
\end{cases}
\]
On $L'/L\cong \F_p$ the quadratic form is equivalent to $q(x)=\alpha x^2/p$, where $\alpha\in \F_p\setminus\{0\}$ is either a square or a non-square, i.e.~the type of $q$ is determined by $\epsilon=\chi_p(\alpha)$.
We use Milgrams formula for the general Gauss sum
\begin{equation}\label{mil}
\sum_{\gamma\in L'/L}e(q(\gamma))=\sqrt{|L'/L|}e(r/8),
\end{equation}
where $e(\tau)=e^{2\pi i \tau}$ as usual.
The left hand side is equal to the standard Gauss sum $\sum_{x(p)} e(\alpha x^2/p)=\chi_p(\alpha)\eps_p\sqrt{p}$. Inserting this into (\ref{mil}) we find that $\epsilon\eps_p=e(r/8)$. We obtain the following table for $r$ modulo $8$:
\begin{table}[h]
\begin{tabular}{l||r|r}
$p\pmod{4}$   & \quad 1 & \quad 3  \\
\hline\hline
$\epsilon=+1$ & $0$ & $2$  \\ \hline
$\epsilon=-1$ & $4$ & $6$  
\end{tabular}
\end{table}

\bigskip

Let $T=\kzxz{1}{1}{0}{1}$ and $S=\kzxz{0}{-1}{1}{0}$ denote the standard generators of $\Sl_2(\Z)$, and write $(\frake_\gamma)_{\gamma\in L'/L}$ for the standard basis of the group algebra $\C[L'/L]$.
Recall that there is a unitary representation $\rho=\rho_L$ of $\Sl_2(\Z)$ on $\C[L'/L]$ given by
\begin{align*}
\rho(T)\frake_\gamma & = e(q(\gamma)) \frake_\gamma,\\
\rho(S)\frake_\gamma & = \frac{i^{-r/2}}{ \sqrt{|L'/L|}}\sum_{\delta\in L'/L} e(-(\gamma,\delta))\frake_\delta.
\end{align*}
This is essentially the Weil representation corresponding to the quadratic module $(L'/L,q)$. For further properties of $\rho$ we refer to \cite{Sh,Bo1,Br1,Od}.
The negative identity matrix acts as $\rho(-E)\frake_\gamma = (-1)^{r/2}\frake_{-\gamma}$. 

Let $k\in \Z$. Similarly as in \cite{Bo1} and \cite{Br1} we denote by $\calA_{k,\rho}$ the space of {\em nearly holomorphic} modular forms for $\Sl_2(\Z)$ of weight $k$ with representation $\rho$.
These are $\C[L'/L]$-valued holomorphic functions $F(\tau)=\sum_{\gamma\in  L'/L} F_\gamma(\tau)\frake_\gamma$ on the upper complex half-plane $\H$ satisfying the usual transformation law and having a Fourier expansion
\[
F(\tau)=\sum_{\gamma\in L'/L} \sum_{\substack{n\in \Z+q(\gamma)\\n\gg -\infty}} a(\gamma,n)e(n\tau) \frake_\gamma.
\]
Thus $F$ may have a pole at the cusp $\infty$.
The subspace of holomorphic modular forms (resp.~cusp forms) is denoted by $\calM_{k,\rho}$ (resp.~$\calS_{k,\rho}$).

The transformation behavior under $-E$ implies that the components $F_\gamma$ of $F$ satisfy
\[
F_\gamma=\begin{cases}F_{-\gamma},& k\equiv r/2 \pmod{2},\\
-F_{-\gamma},& k\not\equiv r/2 \pmod{2}.
\end{cases}
\]
{\it For the purposes of this paper we have to assume that $k\equiv r/2 \pmod{2}$.}

\section{Scalar valued modular forms}

We normalize the usual weight $k$ Petersson slash operator on functions $f:\H\to \C$ by
\[
(f \mid_k M)(\tau)=(\det M)^{k/2} (c\tau + d)^{-k} f(M\tau)
\]
for $M=\kabcd\in \Gl_2^+(\R)$. We will often omit the subscript $k$ if it is clear from the context.

We write $A_k(p,\chi_p)$ for the space of {\em nearly holomorphic} modular forms of weight $k$ for the group $\Gamma_0(p)$ with character $\chi_p$.
These are holomorphic functions on $\H$ which satisfy the transformation law
$f\mid_k M=\chi_p(d) f$ for all $\kabcd\in\Gamma_0(p)$ and are meromorphic at the cusps. The subspace of holomorphic modular forms (resp.~cusp forms) is denoted by $M_k(p,\chi_p)$ (resp.~$S_k(p,\chi_p)$).
Moreover, for $\epsilon\in \{\pm 1\}$ we define the subspaces
\[
A^\epsilon_k(p,\chi_p)=\left\{\text{$ f=\sum_{n\in \Z} a(n) q^n\in A_k(p,\chi_p)$;\quad $a(n)=0$ if $\chi_p(n)=-\epsilon$}\right\}.
\]
Here $q=e^{2\pi i \tau}$ as usual. A classical Lemma due to  Hecke (cf.~\cite{Ogg} Lemma 6, p.~32) implies that
\[
A_k(p,\chi_p)=A^+_k(p,\chi_p)\oplus A^-_k(p,\chi_p).
\]
If $f=\sum_{n\in \Z} a(n) q^n\in A_k^\epsilon(p,\chi_p)$, then the Fourier polynomial 
\[
\sum_{\substack{n\in \Z\\n<0}} a(n) q^n
\]
is called the {\em principal part} of $f$. Finally, we define the spaces $M_k^\epsilon(p,\chi_p)$ and $S^\epsilon_k(p,\chi_p)$ analogously. 

Let $F\in \calA_{k,\rho}$ and write $F_\gamma=\sum_{n}a(\gamma,n)q^n$ for its components.
It is well known that the restriction of $\rho$ to $\Gamma_0(p)$ acts on $\frake_0$ by multiplication with the character $\chi(M)=\chi_p(d)$ for $M=\kabcd$ (see \cite{Sh}, \cite{Od}). Therefore the component $F_0$ belongs to $A_k(p,\chi_p)$.
Since the Fricke involution 
\[
W_p=\zxz{0}{-1}{p}{0}
\]
acts on $A_k(p,\chi_p)$, the function $F_0\mid W_p$ is also contained in $A_k(p,\chi_p)$. 
By means of the operator $V_p=\kzxz{p}{0}{0}{1}$ 
it can be rewritten as follows:
\begin{align}
\nonumber
F_0\mid W_p &= F_0\mid S \mid V_p\\
\nonumber
&= \frac{i^{-r/2}}{\sqrt{p}} \sum_{\gamma\in L'/L} F_\gamma  \mid V_p\\
&= i^{-r/2} p^{(k-1)/2}\sum_{\gamma\in L'/L} F_\gamma(p\tau).\label{bu1}
\end{align}
We get the following Lemma.

\begin{lem}\label{Ftof}
The assignment $F\mapsto f$, where
\begin{equation*}
f=  \frac{i^{r/2}}{2} p^{(1-k)/2}F_0\mid W_p = \frac{1}{2}\sum_{\gamma\in L'/L} F_\gamma(p\tau),
\end{equation*}
defines an injective  homomorphism $\calA_{k,\rho}\to A^\epsilon_k(p,\chi_p)$. Here $\epsilon=\chi_p(\alpha)$ is given by the quadratic form on $L'/L$.
The function $f$ has the Fourier expansion
\[
f=\frac{1}{2}\sum_{n\in \Z} \sum_{\substack{\gamma\in L'/L\\ pq(\gamma)\equiv n \;(p)}} a(\gamma,n) q^n.
\]
\end{lem}

Conversely, for given  $f=\sum_{n\in \Z}a(n)q^n \in  A_k(p,\chi_p)$ we define a function $H\in \calA_{k,\rho}$ as follows. The $\C[L'/L]$-valued function $\frake_0 (f\mid W_p)$ can be viewed as  a nearly holomorphic modular form with representation $\rho$ for the group $\Gamma_0(p)$.
Thus the induced function
\[
H = \sum_{M\in \Gamma_0(p)\backslash\Sl_2(\Z)} \left(\rho(M)^{-1}\frake_0\right) f\mid W_p\mid M
\]
belongs to $\calA_{k,\rho}$.
We now compute its Fourier expansion. 
A system of representatives for $\Gamma_0(p)\backslash\Sl_2(\Z)$ is given by 
\[
\zxz{1}{0}{0}{1},\qquad ST^j=\zxz{0}{-1}{1}{j}\qquad\text{for $j=0,\dots,p-1$.}
\]
Therefore
\begin{align*}
H&= \frake_0 f\mid W_p
+ \sum_{j\;(p)} \left(\rho(ST^j)^{-1}\frake_0\right) f \mid \zxz{-1}{-j}{0}{-p}\\
&=   \frake_0 f\mid W_p + (-1)^k p^{-k/2} \sum_{j\;(p)}\left(\overline{\rho(T^j)}\,\overline{\rho(S)} \frake_0 \right) f(\tau/p+j/p)\\
&= \frake_0 f\mid W_p + (-1)^k p^{-k/2} \frac{i^{r/2}}{\sqrt{p}}\sum_{\gamma\in L'/L} \frake_\gamma \sum_{j\;(p)}e(-q(\gamma)j)f(\tau/p+j/p).
\end{align*}
Using the assumption $k\equiv r/2\pmod{2}$ and inserting the Fourier expansion of $f$, we find
\begin{align*}
H&= \frake_0 f\mid W_p
+i^{-r/2}p^{-1/2-k/2}\sum_{\gamma\in L'/L} \frake_\gamma \sum_{n\in \Z} a(n) e(n\tau/p)\sum_{j\;(p)}e\big((n-pq(\gamma))j/p\big)\\
&= \frake_0 f\mid W_p
+i^{-r/2}p^{1/2-k/2}\sum_{\gamma\in L'/L} \frake_\gamma \sum_{\substack{n\in \Z\\ n\equiv pq(\gamma)\;(p)}} a(n) e(n\tau/p).
\end{align*}
We obtain the following Proposition.

\begin{prop}\label{premainres}
Let $f=\sum_{n}a(n)q^n\in A_k(p,\chi_p)$. Then the function
\begin{equation}\label{defF}
F=\sum_{\gamma\in L'/L} \frake_\gamma F_\gamma =i^{r/2}p^{k/2-1/2}\sum_{M\in \Gamma_0(p)\backslash\Sl_2(\Z)} \left(\rho(M)^{-1}\frake_0\right) f\mid W_p\mid M.
\end{equation}
belongs to $\calA_{k,\rho}$. The components $F_\gamma$ have the Fourier expansion
\begin{align}
\label{bu2}
F_0 &=  \sum_{\substack{n\in \Z\\ n\equiv 0\;(p)}} a(n) e(n\tau/p)
+ i^{r/2}p^{k/2-1/2} f\mid W_p,\\
\label{bu3}
F_\gamma &=  \sum_{\substack{n\in \Z\\ n\equiv pq(\gamma)\;(p)}} a(n) e(n\tau/p)\qquad(\gamma\neq 0).
\end{align}
\end{prop}

If $f\in A^\epsilon_k(p,\chi_p)$, then the Fourier expansion of $F_0$ can be simplified, avoiding the term $f\mid W_p$ (i.e.~the Fourier expansion of $f$ at the cusp $0$). To this end, similarly as in \cite{Kr}, we first characterize $A^\epsilon_k(p,\chi_p)$ using the Hecke operator $U_p$  defined by
\[
f\mid U_p = \sum_{j \;(p)}  f \mid \zxz{1}{j}{0}{p}.
\]
Observe that
\[
f\mid U_p = p^{1-k/2} \sum_{\substack{n\in \Z\\ n\equiv 0\;(p)}} a(n) e(n\tau/p)
\]
for $f$ as above. So our normalization of $U_p$ 
is slightly different than as usual.

\begin{lem}\label{char}
Let $f=\sum_{n\in \Z}a(n) q^n \in  A_k(p,\chi_p)$ and $\epsilon\in \{\pm 1\}$. Then $f$ belongs to
$A_k^\epsilon(p,\chi_p)$, if and only if
\begin{equation}\label{char1}
f \mid U_p =\epsilon \eps_p \sqrt{p} \,f\mid W_p.
\end{equation}
\end{lem}

\begin{proof}
The function $h= f\mid U_p\mid W_p$ is contained in  $A_k(p,\chi_p)$ and (\ref{char1}) is equivalent to
\[
h =\epsilon \overline{\eps_p} \sqrt{p} \,f.
\]
We have
\begin{align*}
h &= \sum_{j \;(p)}  f \mid \zxz{1}{j}{0}{p}\mid \zxz{0}{-1}{p}{0} \\
&= f\mid W_p \mid V_p + \sum_{j \;(p)^*}  f \mid \zxz{j}{-1}{p}{0}\mid \zxz{p}{0}{0}{1},
\end{align*} 
where the summation in $\sum_{j \;(p)^*}$ runs through all primitive residues modulo $p$.
For a given $j\in \Z$ that is coprime to $p$ let $b,d\in\Z$ such that $jd-pb=1$. Then $\kzxz{j}{b}{p}{d}\in \Gamma_0(p)$ and
\begin{equation}\label{h1}
\zxz{j}{-1}{p}{0} =  \zxz{j}{b}{p}{d} \zxz{1}{-d}{0}{p}.
\end{equation}
Thus
\begin{align*}
h &=  f\mid W_p \mid V_p + \sum_{d \;(p)^*} \chi_p(d) f \mid \zxz{p}{-d}{0}{p}\\
&=f\mid W_p \mid V_p + \sum_{n\in \Z} a(n)q^n \sum_{d\;(p)^*} \chi_p(d) e(-nd/p).
\end{align*} 
If we insert the value of the latter Gauss sum, we obtain
\begin{equation}\label{id1}
h =  f\mid W_p\mid V_p + \overline{\eps_p}\sqrt{p} \sum_{n\in  \Z} \chi_p(n)a(n)q^n.
\end{equation}
From this identity the assertion can be deduced. For the implication ``$\Rightarrow$'' we additionally have to use the fact that a nearly holomorphic modular form $g\in  A_k(p,\chi_p)$ with Fourier coefficients $c(n)$ vanishes identically, if $c(n)=0$ for all $n$ coprime to $p$ (\cite{Ogg} Lemma 6).
\end{proof}

An alternative proof of this lemma can be obtained by considering the vector valued function $F$ attached to $f$ via Proposition \ref{premainres}.
One has to insert (\ref{bu2}) and (\ref{bu3}) into the formula (\ref{bu1}) for $F_0\mid W_p$ and carefully compare Fourier expansions on both sides.

\begin{cor} 
\label{cusps}
Let $f$ be a holomorphic function on $\H$ that has the transformation behavior and the Fourier expansion of an element of $A^\epsilon_k(p,\chi_p)$. Suppose that $f$ is meromorphic (or holomorphic or vanishes, respectively) at the cusp $\infty$. Then it is also meromorphic (or holomorphic or vanishes, respectively) at the cusp $0\sim W_p\infty$.
\end{cor}

\begin{thm}\label{mainres}
Let $f=\sum_{n}a(n)q^n\in A^\epsilon_k(p,\chi_p)$ and define $F$ by (\ref{defF}) as before. Then $F\in \calA_{k,\rho}$ and the components $F_\gamma$ have the Fourier expansion
\begin{align}\label{fex1}
F_0 &=  2\sum_{\substack{n\in \Z\\ n\equiv 0\;(p)}} a(n) e(n\tau/p),\\
\label{fex2}
F_\gamma &=  \sum_{\substack{n\in \Z\\ n\equiv pq(\gamma)\;(p)}} a(n) e(n\tau/p)\qquad(\gamma\neq 0).
\end{align}
The map $f\mapsto F$ and the map described in Lemma \ref{Ftof} are  inverse isomorphisms between $A^\epsilon_k(p,\chi_p)$ and $\calA_{k,\rho}$.
\end{thm}

\begin{proof}
According to Proposition \ref{premainres} we only have to prove (\ref{fex1}).
It suffices to show that
\[
f\mid U_p = i^{r/2}\sqrt{p}\,f\mid W_p.
\]
But this immediately follows from Lemma \ref{char} and the fact that $i^{r/2}=\epsilon\eps_p$.
\end{proof}

\bigskip

From now on we assume that $k\leq 0$. 
For an integer $n$ we define
\begin{equation}\label{sn}
s(n)=\begin{cases} 2,& \text{if $n \equiv 0\pmod{p}$},\\
1,& \text{if $n \not\equiv 0\pmod{p}$}.
\end{cases}
\end{equation}
Moreover, we put $\delta=\chi_p(-1)\epsilon$.
The next theorem gives a criterion for the existence of nearly holomorphic modular forms in $A_k^\epsilon(p,\chi_p)$ with prescribed principal part.
We need to consider the space of modular forms $M_{\kappa}^{\delta}(p,\chi_p)$ of dual weight $\kappa=2-k\geq 2$. 
Recall that there are the $2$ Eisenstein series 
\begin{align}\label{defg}
G_\kappa &=1+\frac{2}{ L(1-\kappa,\chi_p)} \sum_{n=1}^\infty \sum_{d\mid n} d^{\kappa-1}\chi_p(d) q^n,\\
\label{defh}
H_\kappa &=\sum_{n=1}^\infty \sum_{d\mid n} d^{\kappa-1}\chi_p(n/d) q^n
\end{align}
in $M_{\kappa}(p,\chi_p)$ (cf.~\cite{He} Werke p.~818), the former corresponding to the cusp $\infty$, the latter corresponding to the cusp $0$.
The linear combination
\begin{equation}\label{eis}
E_\kappa^\delta = 1+\sum_{n\geq 1} B(n)q^n=  1+\frac{2}{L(1-\kappa,\chi_p)} \sum_{n\geq 1} \sum_{d\mid n} d^{\kappa-1}\left( \chi_p(d) +\delta  \chi_p(n/d)\right) q^n
\end{equation}
belongs to $M_{\kappa}^{\delta}(p,\chi_p)$. 
The space $M_{\kappa}^{\delta}(p,\chi_p)$ 
can be decomposed as a direct sum
\[
M_{\kappa}^{\delta}(p,\chi_p) = \C E_\kappa^\delta \oplus S_{\kappa}^{\delta}(p,\chi_p).
\]

\begin{thm}\label{crit}
There exists a nearly holomorphic modular form $f\in A_k^\epsilon(p,\chi_p)$ with prescribed principal part $\sum_{n<0}a(n)q^n$ (where $a(n)=0$ if $\chi_p(n)=-\epsilon$), if and only if
\[
\sum_{n<0} s(n) a(n)b(-n)=0
\]
for every cusp form $g=\sum_{m>0} b(m)q^m$ in $S_\kappa^{\delta}(p,\chi_p)$.
The
constant term $a(0)$ of $f$ is given by the coefficients of the Eisenstein series $E_\kappa^\delta$:
\[
a(0)=- \frac{1}{2}\sum_{n<0} s(n) a(n)B(-n).
\]
\end{thm}

\begin{proof}
Beside $L$ we consider the lattice $L(-1)$, which is given by $L$ as a $\Z$-module but equipped with the quadratic form $-q(\cdot)$. The type of this quadratic form on $L'/L$ is obviously determined by the sign $\delta=\chi_p(-1)\epsilon$.
Observe that the representation $\rho_{L(-1)}$ attached to $L(-1)$ is equal to the dual representation  $\bar{\rho}=\bar{\rho}_L$ of $\rho_L$.
Thus, by Theorem \ref{mainres} we find that the space $\calM_{\kappa,\bar{\rho}}$ of vector valued modular forms of weight $\kappa$ with representation $\bar{\rho}$ can be identified with $M_\kappa^{\delta}(p,\chi_p)$.
 
On the other hand we already know that $\calA_{k,\rho}$ can be identified with  $A_k^{\epsilon}(p,\chi_p)$.

According to Borcherds' duality theorem (Theorem 4.1 in \cite{Bo2}, see also \cite{Br1} chapter 1.3), the obstructions to finding modular forms in $\calA_{k,\rho}$ with prescribed principal part are given by modular forms in $\calM_{\kappa,\bar{\rho}}$.
If we work out all identifications explicitly, we obtain the stated result.
\end{proof}

If one wants to use Theorems \ref{mainres} and \ref{crit} in the context of Borcherds' theory of automorphic products, then the following two propositions will also be important (see section \ref{sec4}). Notice that from the outset analogous results are {\em not} available for the corresponding spaces of vector valued modular forms (see \cite{Bo2} p.~227).

For $f=\sum a(n)q^n\in A_k(p,\chi_p)$ and a Galois automorphism $\sigma\in \Gal(\C/\Q)$ we define the $\sigma$-conjugate of $f$ by
\[
f^\sigma = \sum_{n\in \Z}a^\sigma(n)q^n.
\]
Here $a^\sigma (n)$ denotes the conjugate of $a(n)$.
It is  well known that the spaces $M_\kappa(p,\chi)$, where $\chi$ denotes a quadratic character, have a basis consisting of modular forms with integral rational coefficients (cf.~\cite{DI} Corollary 12.3.8, Proposition 12.3.11).
This implies that $f^\sigma\in A_k(p,\chi_p)$.

\begin{prop}\label{zbasis}
The space $M_\kappa^\delta(p,\chi_p)$ has a basis of modular forms with integral rational coefficients.
\end{prop}

\begin{proof}
It is easily seen that $M_\kappa^\delta(p,\chi_p)$ has a basis of modular forms with coefficients in the ring of integers of a fixed number field $K$. The assertion follows from the fact that the Galois group of $K/\Q$ acts on $M_\kappa^\delta(p,\chi_p)$.
\end{proof}

\begin{prop}\label{ratcoeff}
Let $f=\sum a(n)q^n\in A_k^\epsilon(p,\chi_p)$ and suppose that $a(n)\in \Q$ for $n<0$.
Then all coefficients $a(n)$ are rational and have bounded denominator (i.e.~there is a positive integer $c$ such that $c f$ has coefficients in $\Z$).
\end{prop}

\begin{proof}
Let $\sigma\in \Gal(\C/\Q)$. Then $h=f-f^\sigma$ lies in $A_k^\epsilon(p,\chi_p)$. The assumption on the coefficients $a(n)$ with $n<0$ implies that $h$ is holomorphic at the cusp $\infty$. 
By Corollary \ref{cusps} it is also holomorphic at the cusp $0$ and therefore contained in $M_k^\epsilon(p,\chi_p)$. Since $k\leq 0$, it has to vanish identically. 

Varying $\sigma$, we find that $f$ is invariant under $\Gal(\C/\Q)$ and therefore has rational coefficients. Since the product $f\Delta^N$ of $f$ with a large power of the delta function $\Delta=q\prod_{n\geq 1}(1-q^n)^{24}$ is contained in $M_{k+12N}(p,\chi_p)$, it has coefficients with bounded denominator. Hence $f$ itself has coefficients with bounded denominator.
\end{proof}

\section{Borcherds products on Hilbert modular surfaces}
\label{sec4}

The results of the previous section can be used to obtain a very smooth and explicit formulation of the Borcherds lifting (Theorem 13.3 in \cite{Bo1}, and \cite{Bo2}) for lattices of odd prime determinant and signature $(2,b^-)$. 
To simplify the presentation and notation we only illustrate this in the special $\Orth(2,2)$-case of Hilbert modular surfaces. 

Let $p\equiv 1\pmod{4}$ be an odd prime and 
$K=\Q(\sqrt{p})$ the real quadratic field of discriminant $p$. Let $\calO$ be the ring of integers and $\frakd=(\sqrt{p})$ the different in $K$. We write $x\mapsto x'$ for the conjugation, $\norm(x)=xx'$ for the norm in $K$, and $\tr(x)=x+x'$ for the trace in $K$.

Then $L=\Z^2\oplus\calO$, with the quadratic form $q(a,b,\lambda)=\norm(\lambda)-ab$ for $(a,b,\lambda)\in L$, is an even lattice of signature $(2,2)$. The dual lattice is $L'=\Z^2\oplus\frakd^{-1}$ and the quadratic form on $L'/L\cong\F_p$ represents the squares, i.e.~$\epsilon=+1$.

By Theorem \ref{mainres} the space $\calA_{0,\rho}$ of modular forms with representation $\rho$ of weight $k=0$ is isomorphic to $A^+_0(p,\chi_p)$, and  the space $\calM_{2,\bar{\rho}}$ of holomorphic modular forms with dual representation $\bar{\rho}$ of weight $\kappa=2$ is isomorphic to $M^+_2(p,\chi_p)$.  

The Hilbert modular group $\Gamma_K=\Sl_2(\calO)$ acts on the product of two upper half planes in the usual way. We use $(z_1,z_2)$ as a standard variable on $\H\times\H$ and write $(y_1,y_2)$  for its imaginary part.
Recall that for every positive integer $m$ the subset
\[
\bigcup_{\substack{(a,b,\lambda)\in L'\\ab-\norm(\lambda)=m/p}}
 \{(z_1,z_2)\in \H\times\H;\quad az_1 z_2 +\lambda z_1+\lambda'z_2+b = 0\}
\]
defines a $\Gamma_K$-invariant algebraic divisor $T(m)$ on $\H\times\H$, the Hirzebruch-Zagier divisor of discriminant $m$. It is the inverse image of an algebraic divisor on the quotient $X_K=(\H\times\H)/\Gamma_K$, which will also be denoted by $T(m)$. Here we understand that all irreducible components of $T(m)$ are assigned the multiplicity $1$. This divisor is non-zero if $\chi_p(m)\neq -1$ and is compact if $m$ is not the norm of an ideal in $\calO$.

The subset
\[
S(m)=\bigcup_{\substack{\lambda\in \frakd^{-1}\\-\norm(\lambda)=m/p}}
\{(z_1,z_2)\in\H\times\H;\quad \lambda y_1 +\lambda'y_2=0\}
\]
of $\H\times\H$ is a union of hyperplanes of real codimension $1$. It is invariant under the stabilizer of the cusp $\infty$. 
For a subset $W\subset\H\times\H$ and $\lambda\in \frakd^{-1}$ we write $(W,\lambda)>0$, if $\lambda y_1+\lambda'y_2>0$ for all $(z_1,z_2)\in W$.

For basic facts on Hilbert modular forms we refer to \cite{Fr}, \cite{Ge}.
It is well known that Hilbert modular forms on $\H\times\H$ for the group $\Gamma_K$ can be identified with modular forms for the orthogonal group of the lattice $L$. (In this identification the Hirzebruch-Zagier divisor $T(m)$ essentially corresponds to the Heegner divisor $y_{-m,\gamma}$ in the terminology of Borcherds \cite{Bo2}.)    

Using Theorem \ref{mainres}, Proposition \ref{ratcoeff}, and the above identifications, we may restate Theorem 13.3 in \cite{Bo1} as follows:
 
\begin{thm}\label{hilbert}
Let $f=\sum_{n\in \Z}a(n)q^n\in A_0^+(p,\chi_p)$ and assume that $s(n)a(n)\in \Z$ for all $n<0$ (where $s(n)$ is defined by (\ref{sn})). Then there is a meromorphic function $\Psi(z_1,z_2)$ on $\H\times\H$ with the following properties:

1. $\Psi$ is a meromorphic modular form for $\Gamma_K$  with some unitary character of finite order. The weight of $\Psi$ is equal to the constant coefficient $a(0)$ of $f$. It can also be computed using Theorem \ref{crit}.

2. The divisor of $\Psi$ is determined by the principal part of $f$. It equals \[\sum_{n<0} s(n)a(n) T(-n).\]

3. Let $W\subset\H\times\H$ be a Weyl chamber attached to $f$, i.e.~a connected component of 
\[
\H\times\H-\bigcup_{\substack{n<0\\ a(n)\neq 0}}S(-n);
\]
and put $N=\min \{n;\; a(n)\neq 0\}$.
The function $\Psi$ has the Borcherds product expansion
\[
\Psi(z_1,z_2)=e(\rho_W z_1 + \rho'_W z_2) \prod_{\substack{\nu\in\frakd^{-1} \\ (\nu,W)>0}} \left(1-e(\nu z_1 +\nu' z_2)\right)^{s(p\nu\nu')a(p\nu\nu')}.
\]
Here $\rho_W$ and  $\rho_W'$ are algebraic numbers in $K$ that can be computed explicitly.
The product converges normally for all $(z_1,z_2)$ with $y_1 y_2 > |N|/p$ outside the set of poles. 

4. There exists a positive integer $c$ such that $\Psi^c$ has integral rational Fourier coefficients with greatest common divisor $1$.

\end{thm}

Note that automorphic products for $\Gamma_K$ are constructed in a different way in \cite{Br2}. There the exponents are given by the coefficients of nearly holomorphic Poincar\'e series of weight $2$. Using non-holomorphic Poincar\'e series of weight $0$ as in \cite{Br1}, these coefficients could be related to the coefficients of nearly holomorphic modular forms of weight $0$. But this requires a considerable amount of work. 
An important point of the above formulation of Borcherds' theorem is that we also get the last assertion on the integrality properties of the Fourier coefficients of the lifting. This will be vital for some arithmetic applications.

For the computation of the Weyl vector $(\rho_W,\rho'_W)$ one cannot apply the result of \cite{Bo1} (Theorem 10.4), since the lattice $\calO$ is not isotropic.
However, we may use the formula given in \cite{Br2} p.~72 or \cite{Br1} chapter 2.3. 
Let $W$ be as in the theorem.
It turns out that $\rho_W$, $\rho'_W$ are the uniquely determined numbers in $K$ such that
\begin{equation}\label{wv}
\rho_W y_1 + \rho'_W y_2 = \sum_{n<0}s(n)a(n)\sum_{\substack{\lambda\in \frakd^{-1}\\ \lambda >0\\ \norm(\lambda)=n/p}} \min(|\lambda y_1|, |\lambda'y_2|)
\end{equation}
for all $(z_1,z_2)\in W$.
Let $\eps_0>1$ be the fundamental unit of $K$. It has norm $-1$, because $K$ has prime discriminant. For every negative integer $n$ with $a(n)\neq 0$ there are only finitely many $\lambda\in \frakd^{-1}$ such that
$\lambda>0$, $\norm(\lambda)=n/p$, and
\[
\lambda y_1 + \lambda' y_2<0, \qquad \eps_0^2 \lambda y_1 + {\eps_0'}^2\lambda' y_2>0,
\]
for all $(z_1,z_2)\in W$.
Denote the set of these $\lambda$ by $R(W,n)$. By Dirichlets unit theorem we have
\begin{align}
\nonumber
\sum_{\substack{\lambda\in \frakd^{-1}\\ \lambda >0\\ \norm(\lambda)=n/p}} \min(|\lambda y_1|, |\lambda'y_2|)
&=\sum_{\lambda\in R(W,n)}\left( \sum_{n\geq 0} \eps_0^{-2n} \lambda y_1 - \sum_{n> 0} \eps_0^{-2n} \lambda' y_2\right)\\
\nonumber
&=\sum_{\lambda\in R(W,n)}\left( \frac{\eps_0^2 \lambda}{\eps_0^{2}-1}y_1 -\frac{\lambda'}{\eps_0^{2}-1} y_2 \right)\\
\label{wv2}
&=\frac{1}{\tr(\eps_0)}\sum_{\lambda\in R(W,n)}\left( \eps_0 \lambda y_1 + \eps_0'\lambda' y_2 \right).
\end{align}
Inserting this into (\ref{wv}) we obtain a formula for $\rho_W$ and $\rho_W'$. In particular we find that these numbers are conjugate and contained in $(\tr\eps_0)^{-1}\frakd^{-1}$. 

Finally, we remark that by Theorems 8 and 9 of \cite{Br2} any meromorphic modular form for $\Gamma_K$, whose divisor is a linear combination of Hirzebruch-Zagier divisors, is given by a Borcherds product as in the above theorem.

\bigskip

By a classical result of Hecke \cite{He}, the dimension of the space $S_2(p,\chi_p)$ (where $p\equiv 1\pmod{4}$ is a prime) is equal to $2\left[\frac{p-5}{24}\right]$. It is easily seen that the dimension of the obstruction space $S_2^+(p,\chi_p)$ is half the dimension of $S_2(p,\chi_p)$. Hence, $S_2^+(p,\chi_p)=0$ if and only if  $p=5$, $13$, or $17$.
Let us assume that $p$ is one of these primes. Then by
Theorem \ref{crit}  any Fourier polynomial can be realized as the principal part of a nearly holomorphic modular form in $A^+_0(p,\chi_p)$.
If $m$ is  a positive integer with $\chi_p(m)\neq -1$, we write 
\[
f_m=\sum_{n\geq -m}a_m(n)q^n
\]
for the unique element of $A^+_0(p,\chi_p)$, whose principal part is equal to  $s(m)^{-1} q^{-m}$.
The Borcherds lift $\Psi_m$ of $f_m$ is a holomorphic modular form for $\Gamma_K$ of weight
$a_m(0)=-B(m)/2$, where $B(m)$ is the $m$-th coefficient of the Eisenstein series $E^+_2\in M_2^+(p,\chi_p)$.
The divisor of $\Psi_m$ equals $T(m)$.
The functions $f_m$ can be easily constructed. 
We now indicate this in the case $p=5$. 

The normalized Eisenstein series of weight $2$ for $\Gamma_0(5)$ with trivial character is given by
\[
\frakE_2=1+6\sum_{n\geq 1}(\sigma(n)-5\sigma(n/5)) q^n,
\]
where $\sigma(n)=\sum_{d|n}d$ denotes the sum of divisors of $n$.
The Eisenstein series $G_2$ and $H_2$ for $\Gamma_0(5)$ defined by (\ref{defg}) and (\ref{defh}) can be expressed in terms of the eta function $\eta=q^{1/24}\prod_{n\geq 1}(1-q^n)$ as follows:
\[
G_2(\tau)=\eta(\tau)^5/\eta(5\tau), \qquad H_2(\tau)=\eta(5\tau)^5/\eta(\tau).\]
In particular these Eisenstein series do not vanish on $\H$.
Thus
$\frakE_2/H_2\in A_0(p,\chi_p)$ has a first order pole at $\infty$ and is holomorphic on $\H$ and at the cusp $0$. This implies that it equals $f_1$.
Similarly it can be seen that 
\[
f_4=G_2/H_2 (f_1^3+108 f_1) -9 f_1^3 +1128 f_1.
\]  
The function $f_5$ can be constructed as the product $\frac{1}{2}E^+_2(\tau)J(5\tau)$, where $J(\tau)$ is the unique modular form of weight $-2$ for $\Sl_2(\Z)$ that is holomorphic on $\H$ and whose Fourier expansion starts with $q^{-1}+O(1)$.
Now the other functions $f_m$ can be obtained inductively by multiplying the above functions with powers of $j(5\tau)$ and subtracting suitable multiples of the $f_{m'}$ with smaller index $m'$. Here $j=q^{-1}+744+196884q+\dots$ denotes the usual $j$-function. One finds that the first $f_m$ are:
\begin{align*}
f_1&=
q^{-1} + 5 + 11\,q - 54\,q^{4} + 55\,q^{5} + 44\,q^{6} - 395\,q^{9} + 340\,q^{10} + 296\,q^{11} - 1836\,q^{14}+\dots,\\
f_4&=
q^{-4} + 15 - 216\,q + 4959\,q^{4} + 22040\,q^{5}
 - 90984\,q^{6} + 409944\,q^{9} + 1388520\,q^{10} +\dots,\\
f_5&=\tfrac{1}{2} \,q^{-5} + 15 + 275\,
q + 27550\,q^{4} + 43893\,q^{5} + 255300\,q^{6} + 4173825\,q^{9}
 + \dots,\\
f_6&=
q^{-6} + 10 + 264\,q - 136476\,q^{4} + 306360\,q^{
5} + 616220\,q^{6} - 35408776\,q^{9} + \dots,\\
f_9&=
q^{-9} + 35 - 3555\,q + 922374\,q^{4} + 7512885\,q
^{5} - 53113164\,q^{6} + 953960075\,q^{9} +  \dots,\\
f_{10}&= \tfrac{1}{2}\,q^{-10} + 10 + 
3400\,q + 3471300\,q^{4} + 9614200\,q^{5} + 91620925\,q^{6} + 
5391558200\,q^{9} + \dots.
\end{align*}

Since the divisor of the function $\Psi_1$ is equal to $T(1)$, it has to be equal to a multiple of the classical 
Hilbert modular form $\Theta$ constructed by Gundlach in \cite{Gu} as a product of $10$ theta functions. Because $\Psi_1$ has integral coprime  Fourier coefficients, and $\Theta$ has integral Fourier coefficients with greatest common divisor $64$, we find $\Psi_1=\frac{1}{64}\Theta$.
Let $W$ be the Weyl chamber attached to $f_1$ that contains the point $(-i\eps_0',i\eps_0)$, where $\eps_0=\frac{1}{2}(1+\sqrt{5})$ denotes the fundamental unit of $\Q(\sqrt{5})$. Since the set of $\lambda\in \frakd^{-1}$ with norm $-1/5$ is given by $\{\pm \eps_0^{2n}/\sqrt{5};\; n\in \Z\}$, we obtain that
\[
R(W,-1)=\{1/\sqrt{5}\}.
\]
According to (\ref{wv2}) we have $\rho_W=\eps_0/\sqrt{5}$. Thus $\Psi_1$ has the product expansion
\begin{equation}\label{gund}
\Psi_1(z_1,z_2)=e\left( \eps_0 z_1/\sqrt{5} - \eps_0' z_2/\sqrt{5}\right) \prod_{\substack{\nu\in\frakd^{-1} \\ \eps_0\nu'-\eps_0'\nu>0}} \left(1-e(\nu z_1 +\nu' z_2)\right)^{s(5\nu\nu')a_1(5\nu\nu')}.
\end{equation}

If the divisor $T(m)$ on $X_K$ is compact, then $S(m)$ is empty and $\H\times\H$ is the only Weyl chamber for $f_m$. Thus
\begin{equation}\label{comp}
\Psi_m(z_1,z_2)= \prod_{\substack{\nu\in\frakd^{-1} \\ \nu\gg 0}} \left(1-e(\nu z_1 +\nu' z_2)\right)^{s(5\nu\nu')a_m(5\nu\nu')}.
\end{equation}

\end{document}